\documentclass[11pt]{amsart}
\usepackage{amsmath}
\usepackage{amssymb}
\usepackage{amsthm}
\usepackage[margin=1in]{geometry}
\usepackage[colorlinks=true,citecolor=blue,linkcolor=blue,urlcolor=blue]{hyperref}
\usepackage[T1]{fontenc}
\usepackage{lmodern}
\usepackage{microtype}
\usepackage{needspace}
\usepackage{comment}

\newtheorem{theorem}{Theorem}[section]
\newtheorem{lemma}[theorem]{Lemma}
\newtheorem{proposition}[theorem]{Proposition}
\newtheorem{conjecture}[theorem]{Conjecture}

\theoremstyle{definition}
\newtheorem{definition}[theorem]{Definition}

\newtheorem{corollary}[theorem]{Corollary}

\newcommand{\R}{\mathbb R}

\newcommand{\aff}{\operatorname{aff}}
\newcommand{\supp}{\operatorname{supp}}
\newcommand{\dist}{\operatorname{dist}}
\newcommand{\dimH}{\dim_{\mathrm H}}

\DeclareMathOperator{\dir}{dir}

\title{Pinned Simplex Volumes from Thin Sets in Euclidean Spaces}
\author{T. Borges, B. Foster, A. Iosevich, Y. Ou, E. Palsson}

\address[T. Borges]{Mathematisches Institut, University of Bonn,
Endenicher Allee 60, 53115 Bonn, Germany}
\email{tborges@math.uni-bonn.de}

\address[B. Foster]{Department of Mathematics,
Washington University in St. Louis,
St. Louis, MO 63130, USA}
\email{bfoster@wustl.edu}

\address[A. Iosevich]{Department of Mathematics,
University of Rochester, Rochester, NY 14627, USA}
\email{alex.iosevich@rochester.edu}

\address[Y. Ou]{Department of Mathematics,
University of Pennsylvania, Philadelphia, PA 19104, USA}
\email{yumengou@sas.upenn.edu}

\address[E. Palsson]{Department of Mathematics,
Virginia Tech, Blacksburg, VA 24061, USA}
\email{palsson@vt.edu}

\date{}
\begin{document}

\begin{abstract}
Let $E\subset\mathbb R^d$ be a Borel set, with $d\geq2$.
We prove that if $\dim_{\mathrm H}(E)>d/2$, then there exist
distinct points $a,b\in E$ such that the set of triangle areas
with fixed base $a,b$ and third vertex varying in $E$ has positive
Lebesgue measure. More generally, for $2\leq r\leq d$, if
\[
\dim_{\mathrm H}(E)>
\max\left\{r-1,\frac{d-r+2}{2}\right\},
\]
then there are $r$ affinely independent points in $E$ for which the
corresponding base-pinned $r$-simplex volume set has positive
Lebesgue measure.

The geometric input is a thin-plane estimate of Bright, Ortiz,
and Zakharov \cite{BrightOrtizZakharov}. For a suitable $(r-2)$-dimensional anchor spanned by
$E$, this estimate yields simultaneous Frostman bounds for a
projected vertex measure in the normal space of the anchor and
for the distribution of directions determining bases that share
the anchor. The Peres--Schlag machinery \cite{PeresSchlag}, applied to a local family
of distance-to-line maps in this $(d-r+2)$-dimensional normal
space, then produces the desired base.

Finally, we establish an exact correspondence between fixed-base
simplex volumes and pinned distances in the normal space of the
base, prove a converse for matched bases in layered lifts, and
obtain the lower obstruction $\max\{r-1,d/(2r)\}$ for any
universal dimension threshold.
\end{abstract}

\maketitle



\section{Introduction}
\label{sec:introduction}

Let $d\geq2$ and $1\leq r\leq d$. For points
$x_0,\ldots,x_r\in\mathbb R^d$, write
$\operatorname{Vol}_r(x_0,\ldots,x_r)$ for the $r$-dimensional
Euclidean volume of the simplex with those vertices. When $r=1$,
this is the length of the segment joining $x_0$ and $x_1$, while when
$r=2$, it is ordinary triangle area. For a Borel set
$E\subset\mathbb R^d$, define the set of $r$-volumes in $E$ as
\[
V_r(E):=
\{\operatorname{Vol}_r(x_0,\ldots,x_r):x_0,\ldots,x_r\in E\}\subset \R.
\]
The volume map is continuous, so $V_r(E)$ is analytic and hence
Lebesgue measurable; see \cite[Sections 14 and 21]{Kechris}.

The case $r=1$ is the classical distance-set problem. For
$E\subset\mathbb R^d$ and $x\in\mathbb R^d$, write
\[
\Delta(E)=\{|u-v|:u,v\in E\},\qquad
\Delta_x(E)=\{|u-x|:u\in E\}.
\]
Thus $V_1(E)=\Delta(E)$. Falconer's conjecture asserts that if
$\dim_{\mathrm H}(E)>d/2$ then
$\mathcal L^1(\Delta(E))>0$ \cite{Falconer}.
The corresponding pinned conjecture asks whether the same
hypothesis guarantees a point $x\in E$ with
$\mathcal L^1(\Delta_x(E))>0$.
The best known sufficient thresholds for the pinned problem are
$5/4$ in the plane, due to Guth, Iosevich, Ou, and Wang
\cite{GuthIosevichOuWang}, and
$\frac d2+\frac14-\frac{1}{8d+4}$ for $d\geq3$, due to
Du, Ou, Ren, and Zhang \cite{DuOuRenZhang}.
For the unpinned problem, OpenAI recently announced a proof
of the long-standing Falconer conjecture in every dimension
$d\geq2$ \cite{OpenAIFalconer}.

For $r\geq2$, we ask whether one can fix an affinely independent
base $x_0,\ldots,x_{r-1}\in E$ and obtain a positive-measure set
of volumes as the last vertex varies in $E$. This base-pinned
conclusion is stronger than positivity of $\mathcal L^1(V_r(E))$.
Our main result gives the following sufficient condition.

\begin{theorem}[Volumes of base-pinned $r$-simplices]
\label{thm:main}
Let $2\leq r\leq d$, and let $E\subset\mathbb R^d$ be a Borel set satisfying
\[
\dim_{\mathrm H}(E)>
\max\left\{
r-1,
\frac{d-r+2}{2}
\right\}.
\]
Then there exist affinely independent points
$x_0,\ldots,x_{r-1}\in E$ such that
\[
\mathcal L^1
\left(
\{\operatorname{Vol}_r(x_0,\ldots,x_{r-1},x):x\in E\}
\right)>0.
\]
Consequently, $\mathcal L^1(V_r(E))>0.$
\end{theorem}

\Needspace{11\baselineskip}

In the case $r=2$, the Theorem above gives the triangle-area threshold $d/2$ in every
ambient dimension.

\begin{corollary}[Areas of base-pinned triangles]
\label{cor:triangles}
Let $d\geq2$, and let $E\subset\mathbb R^d$ be a Borel set satisfying
\[
\dim_{\mathrm H}(E)>\frac d2.
\]
Then there exist distinct points $a,b\in E$ such that
\[
\mathcal L^1
\left(
\{\operatorname{Vol}_2(a,b,x):x\in E\}
\right)>0.
\]
Consequently, $\mathcal L^1(V_2(E))>0.$
\end{corollary}

Abundance of good pins is useful when assembling larger
configurations with shared vertices. In
\cite{BorgesEtAlGraphs}, abundance and $L^2$-density results
for pinned $r$-stars are used to build $r$-admissible
distance graphs; related ideas enter the study of volume
vectors determined by hypergraphs in \cite{BorgesEtAl}.
A refinement of the proof of Theorem~\ref{thm:main} gives
the following almost-everywhere conclusion for suitably
chosen measures on $E$. When $d=r$, it strengthens the
base-pinned conclusion of Shmerkin and Yavicoli
\cite{ShmerkinYavicoli} while retaining the sharp
threshold $d-1$.

\begin{corollary}[Abundance of good bases]
\label{cor:almost-every-base}
Let $E\subset\mathbb R^d$ be Borel, let $2\leq r\leq d$,
and suppose that
\[
\dim_{\mathrm H}E>
\max\left\{r-1,\frac{d-r+2}{2}\right\}.
\]
For every $\max\left\{r-1,\frac{d-r+2}{2}\right\}
<s<\dim_{\mathrm H}E,$
there exist compactly supported $s$-Frostman probability
measures $\mu_0,\ldots,\mu_{r-1}$, with supports contained
in $E$ and in good position as in Definition \ref{def:good-position}, such that
\[
\mathcal L^1
\left(
\{\operatorname{Vol}_r(a_0,\ldots,a_{r-1},x):x\in E\}
\right)>0
\]
for $(\mu_0\times\cdots\times\mu_{r-1})$-almost every
$(a_0,\ldots,a_{r-1})$.
\end{corollary}

As a simple application, two $r$-simplices sharing a
fixed good base determine a set of volume pairs of
positive $\mathcal L^2$-measure: the two remaining
vertices vary independently, so this set is the
Cartesian product of the corresponding base-pinned
volume set with itself. More complicated configurations may require additional
control of the pins or stability under restrictions of
the vertex measures; see \cite[Appendix~A]{BorgesEtAl}
for related gluing arguments for chains of triangle
areas. We leave these extensions to the interested reader.

The term $r-1$ in Theorem~\ref{thm:main} is unavoidable: a set
contained in an affine $(r-1)$-plane determines only the value zero for its
$r$-dimensional volumes. The threshold in the theorem equals
$r-1$ precisely when $d\leq3r-4$, so it is sharp in that range. When $r=d\geq 2$ ($d$-volumes in $\R^d$) Theorem \ref{thm:main} recovers a known result of Shmerkin and Yavicoli \cite{ShmerkinYavicoli} and for  $r=2$ and $d\geq 3$ (triangles in $\R^d$,$d\geq 3$) it drops the best known threshold from $\frac{d+1}{2}$ \cite{PalssonPsaromiligkos} to $\frac{d}{2}$. We discuss previous work in more detail below.

As observed by Shmerkin and Yavicoli \cite{ShmerkinYavicoli},
the volume of a simplex with a fixed base is a positive constant
multiple of the distance from its remaining unpinned vertex to the affine
span of the base.
For each fixed $r\geq3$ and sufficiently large $d$, our threshold is
\[
\frac{d-r+2}{2}=\frac d2-\frac{r-2}{2}<\frac d2.
\]
Thus distances from $E$ to a suitably chosen affine $(r-1)$-plane spanned by vertices in $E$ can
produce a positive-measure set under weaker dimensional
assumptions than those conjectured to be necessary for distances
from $E$ to a pinned point in it. The lattice construction in
Section~\ref{sec:falconer-connection} also reflects this
distinction: it gives the lower obstruction
$\max\{r-1,d/(2r)\}$ for simplex volumes, and the term $d/(2r)$ is less restrictive than $d/2$, for increasing $r$.

\subsection{Earlier Results}
Early results on simplex volumes considered configurations with
one vertex fixed at the origin. Erdo\u{g}an, Hart, and Iosevich
\cite{ErdoganHartIosevich} obtained a positive-measure set of
triangle areas in $\mathbb R^2$ under the hypothesis
$\dim_{\mathrm H}(E)>3/2$.
Greenleaf, Iosevich, and Mourgoglou
\cite{GreenleafIosevichMourgoglou} established the corresponding
result for tetrahedral volumes in $\mathbb R^3$ when
$\dim_{\mathrm H}(E)>13/5$.
For $d$-simplex volumes with all vertices varying in $E$,
Grafakos, Greenleaf, Iosevich, and Palsson
\cite{GrafakosGreenleafIosevichPalsson} proved positivity
of Lebesgue measure when
$\dim_{\mathrm H}(E)>d-1+\frac{1}{2d}$ for even $d\geq2$,
or $\dim_{\mathrm H}(E)>d-1+\frac{1}{2(d-1)}$ for odd
$d\geq3$, and conjectured that the necessary threshold
$d-1$ should suffice.

The sharp threshold for positive measure in the matching
ambient dimension was obtained by Shmerkin and Yavicoli
\cite{ShmerkinYavicoli}, whose proof also gives a base-pinned
conclusion. In our notation, if $E\subset\mathbb R^d$ is Borel
and $\dim_{\mathrm H}(E)>d-1$, then there are affinely independent
$x_0,\ldots,x_{d-1}\in E$ such that
\[
\{\operatorname{Vol}_d(x_0,\ldots,x_{d-1},x):x\in E\}
\]
has positive Lebesgue measure. The threshold is sharp because
sets contained in affine hyperplanes determine only zero
$d$-dimensional volumes.

Theorem~1.1 of \cite{ShmerkinYavicoli} was stated for
$r$-simplex volumes in $\mathbb R^d$ throughout the range
$2\leq r\leq d$, with the sufficient condition
$\dim_{\mathrm H}(E)>r-1$. Gait\'an Montejo and Palsson
pointed out, and Borges, Foster, Ou, Palsson, and Romero Acosta
subsequently discussed, that the published proof does not
establish the claimed conclusion when the affine span of the
fixed base has codimension greater than one; see
\cite{GaitanMontejoPalsson,BorgesEtAl}. The proof produces many
distances to a hyperplane containing the base, but distance to
that hyperplane need not control distance to the lower-dimensional
affine span of the base. We use the Shmerkin--Yavicoli theorem
only in the matching case $d=r$, where this obstruction is absent.

Gait\'an Montejo and Palsson resolved the low-codimension range
and obtained a nontrivial theorem in every ambient dimension.
Their sufficient condition is
\[
\dim_{\mathrm H}(E)>\max\{r-1,d-r+1\};
\]
see \cite[Theorem~1.1]{GaitanMontejoPalsson}. This is sharp when
$d\leq2r-2$, since the maximum then equals $r-1$. Our threshold
agrees with theirs in that range and is strictly smaller for
$d\geq2r-1$. In particular, Theorem~\ref{thm:main} extends the
sharp range to $d\leq3r-4$.

For doubly pinned triangle areas, Borges, Foster, Ou, Palsson,
and Romero Acosta obtained an abundance conclusion under
\[
\dim_{\mathrm H}(E)>\frac{d+2}{2},\qquad d\geq 3;
\]
see \cite[Corollary~1.3]{BorgesEtAl}. Their theorem gives a
positive-measure area set for almost every pair of pins in
suitably separated positive-measure pieces of a Frostman measure.
They also proved the corresponding abundance statement in the
plane under $\dimH E>1$, which is used in their triangle-chain
applications. The higher-dimensional threshold comes from the two-star
estimates of Iosevich, Pham, Pham, and Shen \cite{IPPS}.
Combining the improved two-star estimates of Borges, Ou, and
Pasquariello \cite{BorgesOuPasquariello} with
\cite[Theorem~3.5]{BorgesEtAl} lowers the sufficient threshold
for the same abundance conclusion to
\[
\alpha_+(d,2):=\frac{d^2+2d+2}{2d+1}<\frac{d+2}{2},
\qquad d\geq3.
\]

The same argument extends to face-pinned $r$-simplex
volumes. Using the $r$-star estimates of
\cite{BorgesOuPasquariello} and the distance-to-volume
transfer mechanism discussed in \cite[Section~3.2]{BorgesEtAl}
gives the corresponding abundance conclusion whenever
$\dim_{\mathrm H}(E)>\alpha_+(d,r)
:=\frac{d^2+dr+r}{2d+1}$, for $3\leq r<d$.

More recently, Palsson and Psaromiligkos
\cite[Theorem~1.7]{PalssonPsaromiligkos} proved that if
$\dim_{\mathrm H}(E)>(d+1)/2$, then any prescribed first vertex
yields a positive-measure triangle-area set for almost every
second vertex with respect to a suitable finite-energy measure
on $E$. For $d\geq3$, this gives the previously best known
sufficient threshold for the existence of a base-pinned
triangle-area set of positive measure.  In dimension $d=3$, the former threshold
equals $2=d-1$; the threshold $d-1$ for base-pinned triangle
areas was also obtained in
\cite[Theorem~1.1]{GaitanMontejoPalsson} and
\cite[Corollary~4.4]{BorgesEtAl}.

For $d=2$, Corollary~\ref{cor:triangles} recovers the
Shmerkin--Yavicoli theorem. For $d\geq3$, it lowers the sufficient
threshold for the existence of a good triangle base from
$(d+1)/2$ to $d/2$.
Several of the preceding results provide stronger control
over the prescribed pins or establish abundance under higher
dimensional assumptions. Corollary~\ref{cor:almost-every-base}
gives an almost-everywhere conclusion at our threshold
for suitably chosen Frostman measures on $E$.

A stronger conclusion asks whether the volume set has
nonempty interior.
Greenleaf, Iosevich, and Taylor
\cite[Theorem~1.2]{GreenleafIosevichTaylor2022} obtained the
threshold $d-1+1/d$ for $d$-simplex volumes with any prescribed
first vertex, for $d\geq3$.
For unpinned triangle areas, the triangle-configuration
theorem of Palsson and Romero Acosta
\cite{PalssonRomeroAcosta} also yields nonempty interior
under $\dim_{\mathrm H}(E)>(2d+3)/3$, for $d\geq4$,
via the Jacobian transfer argument in
\cite[Section~3.3]{BorgesEtAl}.
In the plane, Greenleaf, Iosevich, and Taylor
\cite[Theorem~1.1(i)]{GreenleafIosevichTaylor2024} obtained
interior for pairs of triangle areas at threshold $3/2$,
and hence also for the scalar area set.
More recently, Palsson and Psaromiligkos
\cite[Corollary~1.10]{PalssonPsaromiligkos} established
nonempty interior for $r$-simplex volumes with any prescribed
first vertex under $\dim_{\mathrm H}(E)>(d+r-1)/2$,
for $2\leq r<d$.

Another line of work studies several volumes simultaneously.
McDonald \cite{McDonald} obtained positive-measure results
for planar area-types, which record the signed areas of
parallelograms spanned by all pairs of points in a configuration.
Galo and McDonald \cite{GaloMcDonald} extended this framework
to volume-types in $\mathbb R^d$, recording the determinants
of all $d$-tuples in a configuration and measuring the resulting
sets in their natural parameter spaces.

\subsection{Discussion of our method}

The starting point is the base-times-height identity
\[
\operatorname{Vol}_r(x_0,\ldots,x_r)
=
\frac{1}{r}\operatorname{dist}(x_r,F)
\operatorname{Vol}_{r-1}(x_0,\ldots,x_{r-1}),
\]
where $x_0,\ldots,x_{r-1}$ are affinely independent and
$F$ is their affine span. Thus it suffices to find such
a base in $E$ for which
\[
\{\operatorname{dist}(x,F):x\in E\}
\]
has positive Lebesgue measure.

The work of Bright, Ortiz, and Zakharov
\cite[Theorem~4.1]{BrightOrtizZakharov} provides dimensional lower bounds
for families of affine planes spanned by a subset of $\R^d$.
This suggests applying the Peres--Schlag machinery to
distance maps parametrized by the $(r-1)$-planes spanned
by $E$. To obtain our threshold, however, we work with a
smaller family of planes sharing a fixed
$(r-2)$-dimensional \emph{anchor}.
For $r=2$, this amounts to fixing a point and varying a
line through it, a setting governed by Ren's radial
projection theorem \cite{Ren}. For $r\geq 3$,
the thin-plane estimates underlying the work of
Bright, Ortiz, and Zakharov provide the corresponding
control for families of base planes containing a common anchor. In the irreducible case, which is the only
case needed here, their argument builds on Ren's
theorem through induction on the dimension of the planes.

The geometric advantage of fixing an anchor $A$ is
that, after orthogonal projection onto its normal space
and translation so that $A$ maps to the origin, every
candidate base plane containing $A$ becomes a line
through the origin. Heights above the original base
planes become distances to these lines. We therefore
work with a fixed projected set and a family of
directions determined by the remaining base vertex.

To implement this approach, we first reduce to an
irreducible Frostman measure with exponent $s$ above
the threshold in Theorem~\ref{thm:main}. If the resulting
ambient dimension equals $r$, we apply the
Shmerkin--Yavicoli theorem. In the remaining case $d>r$, the thin-plane
estimates allow us to choose an anchor $A$ spanned by
$r-1$ vertices of $E$ and such that upon projection to $A^{\perp}$, we have $\sigma$-dimensional families of both points and lines, where 
\[
\frac{d-r+2}{2}<\sigma<\min\{s,d-r+1\}.
\]
This is made precise using certain pushforward measures, and the same thin-plane bound controls both measures:
projected balls lift to neighborhoods of fibers
contained in candidate base planes, while direction
balls lift to neighborhoods of those planes inside
the bounded support.

A local transversality calculation for distance-to-line
maps in the $(d-r+2)$-dimensional normal space, combined
with the Peres--Schlag theorem \cite{PeresSchlag},
bounds the Hausdorff dimension of the directions
for which the distance-to-line image has zero
Lebesgue measure strictly below $\sigma$. Since $\sigma$ is sufficiently large, we can avoid these directions, and so, a suitable base can
be chosen from $E$. This treats all $r\geq2$ uniformly
in the case $d>r$; when $r=2$, the anchor is a single
point and we don't need to project.

Section~\ref{sec:geometric-input} records the thin-plane
input. Section~\ref{sec:common-preliminaries} collects
the analytic preliminaries, and
Section~\ref{sec:transversality} verifies transversality
for distances to lines. The main proof is given in
Section~\ref{sec:main-proof}. Finally,
Section~\ref{sec:falconer-connection} develops the exact
connections with pinned distances and the lattice
obstruction.

\subsection*{Acknowledgments}

T. Borges was partially supported by the European Research
Council (ERC) under the European Union's Horizon Europe
programme (Synergy Grant HALF, grant agreement No.~101224275). B. Foster was partially supported by an AMS-Simons travel grant. A. Iosevich was supported in part by NSF grant DMS-2506858 and a Simons Fellowship in Mathematics in 2026. Y. Ou was supported in part by NSF grant DMS-2142221.

\section{Geometric input}
\label{sec:geometric-input}

We recall the dimensional lower bound of Bright, Ortiz, and
Zakharov \cite{BrightOrtizZakharov} for the family of
$(r-1)$-planes spanned by a set, and then state the
thin-plane estimate used in our proof.

Write $\aff(x_0,\ldots,x_{r-1})$ for the affine span of
$x_0,\ldots,x_{r-1}$. A finite Borel measure $\mu$ is
\emph{irreducible} in an affine space $H$ if it is supported
on $H$ and gives zero mass to every proper affine subspace
of $H$. When $H=\mathbb R^d$, we simply call $\mu$
irreducible. Recall also that $\mu$ is $s$-Frostman if
\[
\mu(B(x,\rho))\leq C\rho^s
\]
for some constant $C$ and all $x\in\mathbb R^d$ and $\rho>0$.

\begin{theorem}[Bright--Ortiz--Zakharov,
{\cite[Theorem~4.1]{BrightOrtizZakharov}}]
\label{thm:boz-dimension}
Let $2\leq r\leq d$ and $0<s\leq d$. If a Borel set
$E\subset\mathbb R^d$ supports an irreducible $s$-Frostman
probability measure, then
\[
\dim_{\mathrm H}
\bigl\{
\aff(x_0,\ldots,x_{r-1}):
x_0,\ldots,x_{r-1}\in E
\text{ are affinely independent}
\bigr\}
\geq r\min\{s,d-r+1\}.
\]
Here the dimension is taken in $A(d,r-1)$, the affine Grassmannian
of $(r-1)$-planes in $\mathbb R^d$.
\end{theorem}

Our proof uses the quantitative thin-plane estimates
underlying this theorem. We first introduce the relevant
geometric notation.

For a linear subspace $W\subset\mathbb R^d$, let $\pi_W$
denote orthogonal projection onto $W$. Let $\mathcal{L}:=G(W,1)$ be
the space of unoriented lines through the origin in $W$,
equipped with the metric
\[
d_{\mathcal{L}}(\ell,\ell')
:=
\|\pi_\ell-\pi_{\ell'}\|_{\mathrm{op}},
\]
where $\|\cdot\|_{\mathrm{op}}$ is the operator norm.
For an affine subspace $A$, write $\dir(A)$ for its
parallel linear subspace and
\[
A(\rho)
:=
\{x\in\mathbb R^d:\dist(x,A)<\rho\},
\qquad \rho>0.
\]
Let $S_{r-1}\subset(\mathbb R^d)^r$ be the set of
affinely dependent $r$-tuples.

\begin{definition}[Good position]
\label{def:good-position}
Compactly supported probability measures
$\mu_0,\ldots,\mu_{r-1}$, with $E_i=\supp\mu_i$, are
in $C$-good position if their supports lie in a common
ball of radius $C$ and
\[
\dist(E_0\times\cdots\times E_{r-1},S_{r-1})
\geq C^{-1},
\]
where distance is measured in the Euclidean product space
$(\mathbb R^d)^r$. Thus every tuple with one vertex in
each support is uniformly affinely independent.
We say that the measures are in good position when
the value of $C$ is not important.
\end{definition}

An irreducible measure on $\mathbb R^d$ has $r$ affinely
independent support points. Restricting and normalizing
it in sufficiently small closed balls around these
points gives measures in good position. These
restrictions remain irreducible and retain the same
Frostman exponent.

\begin{lemma}[Thin-plane input]
\label{lem:boz-thin-input}
Let $2\leq r\leq d$ and $0<s\leq d$. Let
$\mu_0,\ldots,\mu_{r-1}$ be compactly supported,
irreducible $s$-Frostman probability measures in good
position. For every
\[
0<\sigma<\min\{s,d-r+1\}
\quad\text{and}\quad
0<\varepsilon<1,
\]
there are a Borel set
\[
G\subset
\supp\mu_0\times\cdots\times\supp\mu_{r-1}
\]
and a constant $C_\varepsilon$ such that
\[
(\mu_0\times\cdots\times\mu_{r-1})(G)
\geq 1-\varepsilon
\]
and
\begin{equation}
\label{eq:boz-input}
\mu_i\bigl(\aff(x_0,\ldots,x_{r-1})(\rho)\bigr)
\leq C_\varepsilon\rho^\sigma
\qquad
(0\leq i\leq r-1,\ \rho>0)
\end{equation}
for every $(x_0,\ldots,x_{r-1})\in G$.
\end{lemma}

\begin{proof}
After a similarity, all supports lie in the unit ball.
Set
\[
s_0:=\min\{s,d-r+1\}.
\]
The measures are also $s_0$-Frostman, with
$0<\sigma<s_0\leq d-1$.
We apply \cite[Theorem~4.4]{BrightOrtizZakharov}
to $(r-1)$-planes, with Frostman exponent $s_0$,
thin-plane exponent $\sigma$, and error $\varepsilon$.

The quantitative irreducibility hypothesis in that
theorem requires suitable $w,\tau>0$ such that
\[
\mu_i(H(w))\leq\tau
\]
for every proper affine subspace $H$ and every
$0\leq i\leq r-1$. First choose $\tau$ below the
threshold required by the theorem. By
\cite[Lemma~4.7]{BrightOrtizZakharov}, irreducibility
then allows a common $w>0$, small enough for all the
measures and for the theorem. Its conclusion, together
with \cite[Definition~3.3]{BrightOrtizZakharov}, gives
\eqref{eq:boz-input} on a Borel set of product measure
at least $1-\varepsilon$. Undoing the similarity changes
only the constant.
\end{proof}

The estimate controls the full measure $\mu_i$ near each
plane determined by a tuple in $G$. By Fubini, we may
fix $x_0,\ldots,x_{r-2}$ while retaining a positive-mass
set of choices for $x_{r-1}$. The resulting base planes
all contain the anchor
\[
A=\aff(x_0,\ldots,x_{r-2}).
\]
In the main proof, this estimate supplies Frostman bounds
both for the projected vertex measure and for the
distribution of directions extending $A$.
When $r=2$, the anchor is a point and the base planes
are lines. This thin-line case follows from Ren's
discretized radial projection theorem \cite{Ren},
which builds on the radial projection methods of
Orponen, Shmerkin, and Wang \cite{OrponenShmerkinWang},
and forms the starting point of the induction in
\cite{BrightOrtizZakharov}.

\section{Analytic preliminaries}
\label{sec:common-preliminaries}

We state the external results in the forms used below. Standard facts about
Hausdorff measure, Frostman measures, and energy integrals can be found in
\cite{MattilaGeometry}.

\subsection{Frostman measures}

Here and below, $B(x,q)$ denotes the open Euclidean ball centered at $x$ with
radius $q$.

\begin{lemma}
\label{lem:frostman}
Let $F\subset\mathbb R^d$ be compact, and suppose that
$0<t<s<\dim_{\mathrm H}(F)$. Then there is a Borel probability measure
$\mu$ supported on $F$ and a constant $C$ such that
\begin{equation}
\mu(B(x,q))\leq Cq^s
\label{eq:frostman}
\end{equation}
for every $x\in\mathbb R^d$ and every $q>0$.
Moreover,
\begin{equation}
I_t(\mu)
=
\int\int |x-y|^{-t}d\mu(x)d\mu(y)
<\infty,
\label{eq:energy}
\end{equation}
and $\mu(A)=0$ whenever $A\subset\mathbb R^d$ is Borel with
$\dim_{\mathrm H}(A)<s$.
\end{lemma}

\begin{proof}
Frostman's lemma gives a probability measure satisfying
\eqref{eq:frostman}; see \cite[Theorem 8.8]{MattilaGeometry}. Since $s>0$, this measure has no atoms. For a fixed
$x$, decompose the region $|x-y|<1$ into the annuli
\[
2^{-j-1}\leq |x-y|<2^{-j},
\]
where $j=0,1,2,\ldots$. On the $j$th annulus,
$|x-y|^{-t}\leq2^{(j+1)t},$
while
$\mu(B(x,2^{-j}))\leq C2^{-js}.$

It follows that
\[
\int_{|x-y|<1}|x-y|^{-t}d\mu(y)
\leq
C2^t\sum_{j=0}^{\infty}2^{-j(s-t)}.
\]
The part over $|x-y|\geq1$ is at most $1$. The geometric series converges
because $s>t$, and integration in $x$ proves \eqref{eq:energy}.

If $\dim_{\mathrm H}(A)<s$, then $\mathcal H^s(A)=0$. Covering
$A$ by balls and using \eqref{eq:frostman} gives
\[
\mu(A)\leq C\mathcal H^s(A)=0.
\]
\end{proof}

\subsection{The Peres--Schlag exceptional-set theorem}

Let $\Omega$ be a compact subset of $\mathbb R^d$, let
$\Lambda\subset\mathbb R^N$ be open, and let
\[
P_\lambda:\Omega\longrightarrow\mathbb R
\]
be a family that is smooth in $\lambda$. For distinct $x,y\in\Omega$, put
\[
\Phi_\lambda(x,y)
=
\frac{P_\lambda(x)-P_\lambda(y)}{|x-y|}.
\]
The family is transversal of degree zero if, on every compact subset of
$\Lambda$, there are constants $c_0,c_1>0$ such that
\begin{equation}
|\Phi_\lambda(x,y)|\leq c_0
\Longrightarrow
|\nabla_\lambda\Phi_\lambda(x,y)|\geq c_1
\label{eq:transversality}
\end{equation}
for all distinct $x,y\in\Omega$. It is regular of degree zero if, for every
compact $K\subset\Lambda$ and every parameter multi-index $\beta$, there is
a constant $C_{K,\beta}$ such that
\[
|\partial_\lambda^\beta\Phi_\lambda(x,y)|
\leq
C_{K,\beta}
\]
for every $\lambda\in K$ and all distinct $x,y\in\Omega$. This is slightly
stronger than the finite collection of derivative bounds needed in the
theorem, and it is the form verified for the distance-to-line family below.

We also require uniform smoothness: for every compact $K\subset\Lambda$ and
every multi-index $\beta$,
\[
\sup_{\lambda\in K}\sup_{x\in\Omega}
|\partial_\lambda^\beta P_\lambda(x)|<\infty.
\]

\begin{theorem}[Peres--Schlag]
\label{thm:peres-schlag}
Suppose that the scalar-valued $N$-parameter family $P_\lambda$ satisfies the
uniform smoothness condition and is transversal and regular of degree zero.
Let $\mu$ be a finite Borel measure on $\Omega$ such that
\[
I_t(\mu)<\infty
\]
for some $t>1$. Then
\begin{equation}
\dim_{\mathrm H}
\{\lambda\in\Lambda:(P_\lambda)_\#\mu\not\ll\mathcal L^1\}
\leq
N+1-t.
\label{eq:peres-schlag}
\end{equation}
\end{theorem}

This is the degree-zero scalar case of
\cite[Theorem 7.3]{PeresSchlag}. The explicit bound in this notation also
follows from \cite[Theorem 2.2]{Bushling}. Here the pushforward is defined by
\[
(P_\lambda)_\#\mu(C)=\mu(P_\lambda^{-1}(C))
\]
for Borel sets $C\subset\mathbb R$.

If $(P_\lambda)_\#\mu$ is nonzero and absolutely continuous, then
$P_\lambda(\operatorname{supp}(\mu))$ has positive Lebesgue measure. Indeed,
the pushforward is carried by that compact image, and a nonzero absolutely
continuous measure cannot be carried by a Lebesgue null set.

\section{Transversality for distances to lines through one point}
\label{sec:transversality}

The next lemma is the local geometric input. Its two hypotheses exclude the
two degeneracies that occur in the calculation. If $x_0-a$ is parallel to
$u_0$, then $x_0$ lies on the line and the distance function is not smooth. When $d\geq3$, the jet map below is also singular if
$x_0-a$ is perpendicular to $u_0$.

\begin{lemma}
\label{lem:line-transversality}
Let $d\geq2$. Fix $a,x_0\in\mathbb R^d$ with $x_0\ne a$, and fix
$u_0\in S^{d-1}$. Suppose that $x_0-a$ is neither parallel nor perpendicular
to $u_0$. Then there are an open ball $B$ centered at $x_0$ and a spherical
cap $\mathcal{D}$ centered at $u_0$ such that the family
\[
P_{a,u}(x)=\operatorname{dist}(x,a+\mathbb Ru)
\]
for $u\in\mathcal D$ is uniformly smooth, transversal, and regular of degree zero on
every compact subset of the product of the parameter cap and $B$. The cap has local coordinates
$\lambda\in\Lambda\subset\mathbb R^{d-1}$ for which the coordinate map and
its inverse are bi-Lipschitz on the relevant closures.
\end{lemma}

\begin{proof}
Translation and rotation preserve all the relevant distances. We may assume
that $a=0$ and
\[
u_0=e_d=(0,\ldots,0,1).
\]
Write
\[
x=(x',x_d),
\]
where $x'\in\mathbb R^{d-1}$. The assumptions on $x_0$ give
\[
|x_0'|>0\text{ and }
(x_0)_d\ne0.
\]

Directions near $e_d$ can be parametrized by
\[
u(\lambda)
=
\frac{(\lambda,1)}{(1+|\lambda|^2)^{\frac{1}{2}}},\qquad \lambda\in\mathbb R^{d-1}.
\]
Set
\[
P_\lambda(x)
:=\text{dist}(x,\R u(\lambda))
=\left(
|x|^2-(x\cdot u(\lambda))^2
\right)^{\frac{1}{2}}.
\]
At $\lambda=0$, 
\begin{equation}
P_0(x)=|x'|
\label{eq:line-p-zero}
\end{equation}
and direct differentiation gives for $1\leq j\leq d-1$,
\begin{equation}
\left.\frac{\partial P_\lambda}{\partial\lambda_j}(x)\right|_{\lambda=0}
=
-\frac{x_dx_j}{|x'|}.
\label{eq:line-first-derivative}
\end{equation}

Consider the jet map $J_{\lambda}:\R^d\rightarrow\R^d$
\[
J_\lambda(x)
=
\left(
P_\lambda(x),
\nabla_\lambda P_\lambda(x)
\right).
\]
Equations \eqref{eq:line-p-zero} and
\eqref{eq:line-first-derivative} give
\begin{equation}
J_0(x)
=
\left(
|x'|,-\frac{x_dx'}{|x'|}
\right).
\label{eq:line-jet}
\end{equation}

We verify that $D_xJ_0(x)$ is invertible whenever
$x'\ne0$ and $x_d\ne0$. Writing $\omega=\omega(x')=x'/|x'|$, we have
\[
D_xJ_0(x)
=
\begin{pmatrix}
D_{x'}|x'| & \partial_{x_d}|x'|\\[4pt]
D_{x'}(-x_d\omega(x')) &
\partial_{x_d}(-x_d\omega(x'))
\end{pmatrix}
=
\begin{pmatrix}
\omega^T & 0\\[4pt]
-\dfrac{x_d}{|x'|}
\bigl(I_{d-1}-\omega\omega^T\bigr) & -\omega
\end{pmatrix}.
\]
Thus $h=(h',h_d)$ belongs to the kernel of $D_xJ_0(x)$
if and only if
\[
\omega\cdot h'=0,
\qquad
-\frac{x_d}{|x'|}
\bigl(h'-(\omega\cdot h')\omega\bigr)-h_d\omega=0.
\]
Using the first equation, the second simplifies to
\[
-\frac{x_d}{|x'|}h'-h_d\omega=0.
\]
Taking the inner product with $\omega$ gives $h_d=0$,
and then $x_d\ne0$ implies $h'=0$.
Hence $D_xJ_0(x)$ is invertible.


Let $M=D_xJ_0(x_0)$, and let $\kappa>0$ be its smallest singular value. By
continuity, there are a convex open ball $B$ centered at $x_0$ and a bounded
open parameter neighborhood $\Lambda$ centered at $0$ such that
\[
P_\lambda(x)>0
\]
and
\[
\|D_xJ_\lambda(x)-M\|\leq\frac{\kappa}{2}
\]
for every $x\in\overline B$ and every
$\lambda\in\overline\Lambda$. The first inequality keeps $x$ away from the
line $\mathbb Ru(\lambda)$, where the distance function is not smooth.

For $x,y\in\overline B$, the segment joining them lies in $\overline B$.
The fundamental theorem of calculus gives
\[
J_\lambda(x)-J_\lambda(y)
=
M(x-y)
+
\int_0^1
(D_xJ_\lambda(y+t(x-y))-M)(x-y)dt.
\]
It follows that
\begin{equation}
|J_\lambda(x)-J_\lambda(y)|
\geq
\frac{\kappa}{2}|x-y|.
\label{eq:line-jet-separation}
\end{equation}

For distinct $x,y\in\overline B$, division by $|x-y|$ gives
\[
\frac{J_\lambda(x)-J_\lambda(y)}{|x-y|}
=
(\Phi_\lambda(x,y),\nabla_\lambda\Phi_\lambda(x,y)).
\]
If
\[
|\Phi_\lambda(x,y)|\leq\frac{\kappa}{4},
\]
then \eqref{eq:line-jet-separation} implies
\[
|\nabla_\lambda\Phi_\lambda(x,y)|
\geq
\frac{\kappa}{4}.
\]
This is degree-zero transversality.

It remains to verify regularity. The function
$(x,\lambda)\mapsto P_\lambda(x)$ is smooth on a neighborhood of the compact
set $\overline B\times\overline\Lambda$. For every parameter multi-index
$\beta$, the $x$-gradient of
\[
x\longmapsto\partial_\lambda^\beta P_\lambda(x)
\]
is uniformly bounded there. The mean value theorem gives
\[
|\partial_\lambda^\beta P_\lambda(x)
-\partial_\lambda^\beta P_\lambda(y)|
\leq
C_\beta|x-y|.
\]
After division by $|x-y|$, this supplies all the required regularity bounds
for $\Phi_\lambda$. Smoothness on the same compact set gives uniform
smoothness for $P_\lambda$.

Finally, $\lambda\mapsto u(\lambda)$ is a smooth diffeomorphism from a
sufficiently small neighborhood onto a spherical cap $\mathcal{D}$ about $u_0$. After
shrinking the neighborhoods, the map and its inverse have bounded
derivatives on the relevant closures and are therefore bi-Lipschitz.
\end{proof}

\section{Proof of Theorem~\ref{thm:main}}
\label{sec:main-proof}

We choose one anchor and one restriction of the vertex measure.
The thin-plane estimate gives the two Frostman bounds needed in order to apply
the Peres--Schlag argument directly.

\begin{proof}[Proof of Theorem~\ref{thm:main}]
By compact approximation for Borel sets, we may pass to a compact
subset whose dimension is still above the threshold; see
\cite{Howroyd}. Hence, assume without loss of generality that $E$ is compact. Choose $s$ so that
\[
\max\left\{r-1,\frac{d-r+2}{2}\right\}
<s<\dimH E
\]
and an $s$-Frostman probability measure $\mu$ on $E$.

We may further assume that $\mu$ is irreducible in the
ambient space. Indeed, choose an affine subspace $H$
of minimal dimension with $\mu(H)>0$, and replace $\mu$
by its normalized restriction to $H$. The resulting
measure is irreducible in $H$ and remains $s$-Frostman.
Since $s>r-1$, we have $\dim H\geq r$.
Identifying $H$ isometrically with a Euclidean space
only decreases the ambient dimension. We again denote
that dimension by $d$ and replace $E$ by the support
of the restricted measure. Thus $\mu$ is irreducible
in $\mathbb R^d$ and
\[
s>\max\left\{r-1,\frac{d-r+2}{2}\right\}.
\]

If $d=r$, the conclusion follows from
Shmerkin--Yavicoli~\cite{ShmerkinYavicoli}. Hence, assume $d>r$
and put $d_\perp:=d-r+2\geq3$. Choose $\sigma$ so that
\begin{equation}\label{eq:sigma-choice}
\frac{d_\perp}{2}<\sigma<\min\{s,d_\perp-1\}.
\end{equation}
This is possible because $s>d_\perp/2$ and $d_\perp>2$.

Restricting and normalizing $\mu$ in sufficiently small closed
balls around $r$ affinely independent support points (which exist since $s>r-1$) gives
irreducible $s$-Frostman probabilities
$\mu_0,\ldots,\mu_{r-1}$ in good position. Write
$E_i=\supp\mu_i$. Apply Lemma~\ref{lem:boz-thin-input} with exponent $\sigma$ and $\varepsilon=1/2$.
We obtain a Borel set $G\subset E_0\times\cdots\times E_{r-1}$
of product measure at least $1/2$ such that
\begin{equation}\label{eq:direct-thin}
\mu_{r-1}\bigl(\aff(a_0,\ldots,a_{r-2},x)(\delta)\bigr)
\leq C\delta^\sigma
\end{equation}
for every $(a_0,\ldots,a_{r-2},x)\in G$ and every $\delta>0$.

By Fubini, fix a choice for the first $r-1$ vertices $\mathbf a=(a_0,\ldots,a_{r-2})$ such that the fiber
$G_{\mathbf a}:=\{x\in E_{r-1}:(\mathbf a,x)\in G\}$ has positive
$\mu_{r-1}$-mass. Choose a compact
$K\subset G_{\mathbf a}$ with $\mu_{r-1}(K)>0$, and set
\[
\tilde{\mu}:=\frac{\mu_{r-1}|_K}{\mu_{r-1}(K)},\qquad
A:=\aff(a_0,\ldots,a_{r-2}),\qquad V:=\dir(A).
\]
The fact that the measures $\mu_0,\mu_1,\dots ,\mu_{r-1}$ are in good position implies that the points $a_0,a_1,\dots,a_{r-2},x$  are affinely independent for any $x\in K\subset E_{r-1}$, hence $\dim V=r-2$ and $K\cap A=\varnothing$.

Define the translated orthogonal projection
\[
\Pi_A(x):=\pi_{V^\perp}(x-a_0)
\]
and the line-direction map
\begin{equation}\label{eq: line direction pushforward}
L_{\mathbf a}(x):=\mathbb R\Pi_A(x)\in G(V^\perp,1),
\qquad x\in K.
\end{equation}
Since $K\cap A=\varnothing$, the direction map is
well-defined and continuous on $K$. Define the
pushforward probability measures
\[
\nu:=(\Pi_A)_\#\tilde{\mu}
\quad\text{on }V^\perp,
\qquad
\lambda:=(L_{\mathbf a})_\#\tilde{\mu}
\quad\text{on }G(V^\perp,1).
\]

We claim that both measures are $\sigma$-Frostman.
First let $B(z,\delta)\subset V^\perp$ be a ball whose
preimage meets $K$, and choose
$x_*\in K\cap\Pi_A^{-1}(B(z,\delta))$. If
$y\in\Pi_A^{-1}(B(z,\delta))$, then
\[
\dist(y,x_*+V)=|\Pi_A(y)-\Pi_A(x_*)|<2\delta.
\]
Since $x_*+V\subset\aff(A,x_*)$, it follows that
\[
\Pi_A^{-1}(B(z,\delta))
\subset(x_*+V)(2\delta)
\subset\aff(A,x_*)(2\delta).
\]
By \eqref{eq:direct-thin},
\begin{equation}\label{eq:projected-thin-frostman}
\nu(B(z,\delta))
\leq\frac{\mu_{r-1}(\aff(A,x_*)(2\delta))}{\mu_{r-1}(K)}
\lesssim_{K}\delta^\sigma.
\end{equation}

Next let $\mathcal B(\ell,\delta)$ be a ball in
$G(V^\perp,1)$ whose preimage meets $K$, and choose
\[
x_*\in K\cap
L_{\mathbf a}^{-1}(\mathcal B(\ell,\delta)).
\]
Set $\ell_*:=L_{\mathbf a}(x_*)$ and
$R:=\sup_{y\in K}|\Pi_A(y)|<\infty$.
For every
$y\in K\cap L_{\mathbf a}^{-1}(\mathcal B(\ell,\delta))$,
the triangle inequality gives
\[
d_{\mathcal L}(L_{\mathbf a}(y),\ell_*)
\leq
d_{\mathcal L}(L_{\mathbf a}(y),\ell)
+d_{\mathcal L}(\ell,\ell_*)
<2\delta.
\]
Since $\Pi_A(y)$ lies on $L_{\mathbf a}(y)$, it follows that
\[
\begin{aligned}
\dist(y,\aff(A,x_*))
&=\dist(\Pi_A(y),\ell_*)\\
&=\bigl|(\pi_{L_{\mathbf a}(y)}-\pi_{\ell_*})
          \Pi_A(y)\bigr|\\
&\leq |\Pi_A(y)|\,
d_{\mathcal L}(L_{\mathbf a}(y),\ell_*)\\
&<2R\delta.
\end{aligned}
\]

Applying \eqref{eq:direct-thin} again gives
\begin{equation}\label{eq:direction-thin-frostman}
\lambda(\mathcal B(\ell,\delta))
\leq\frac{\mu_{r-1}(\aff(A,x_*)(2R\delta))}{\mu_{r-1}(K)}
\lesssim_{K}\delta^\sigma.
\end{equation}

If a preimage does not meet $K$, the corresponding
pushforward mass is zero. Thus, both Frostman claims hold. Note that when $r=2$, the Frostman condition for $\lambda$ essentially recovers \cite[Proposition 9.1]{Ren}.

The measure $\nu$ is irreducible in $V^\perp$. Indeed,
$\tilde{\mu}$ is a positive-mass restriction of an irreducible measure,
and the inverse image under $\Pi_A$ of a proper affine subspace
of $V^\perp$ is a proper affine subspace of $\mathbb R^d$.

Choose $\ell_0\in\supp\lambda$ and a unit vector
$u_0\in V^\perp$ such that $\ell_0=\mathbb Ru_0$.
Since $d_\perp\geq3$, both $\mathbb Ru_0$ and $u_0^\perp$
are proper linear subspaces of $V^\perp$, where the
orthogonal complement is taken within $V^\perp$.
Irreducibility therefore allows us to choose
\[
z_0\in\supp\nu\setminus(\mathbb Ru_0\cup u_0^\perp).
\]

Lemma~\ref{lem:line-transversality}, applied in $V^\perp$ with
$a=0$, gives a ball $B\ni z_0$ and a cap $\mathcal D\ni u_0$
on which
\[
P_u(z):=\dist(z,\mathbb R u)
\]
is uniformly smooth, transversal, and regular of degree zero.
Choose a closed ball $B_0\subset B$ centered at $z_0$ with
$\nu(B_0)>0$, and put
\[
\nu_*:=\frac{\nu|_{B_0}}{\nu(B_0)}.
\]
This measure remains $\sigma$-Frostman.

Shrink $\mathcal D$ around $u_0$ so that
$u\mapsto\mathbb Ru$ is bi-Lipschitz onto an open
neighborhood of $\ell_0$ in $G(V^\perp,1)$.
This neighborhood has positive $\lambda$-mass because
$\ell_0\in\supp\lambda$. The restriction of $\lambda$ to
it therefore corresponds to a nonzero $\sigma$-Frostman measure
$\widetilde\lambda$ on $\mathcal D$.

By \eqref{eq:sigma-choice}, choose
\[
d_\perp-\sigma<t<\sigma.
\]
Since $\sigma<d_\perp-1$, we have $t>1$. Moreover,
$I_t(\nu_*)<\infty$ since $t<\sigma$. The cap has $d_\perp-1$ local
parameters, so Theorem~\ref{thm:peres-schlag} gives
\begin{equation}\label{eq:main-exceptional-set}
\dimH\left\{u\in\mathcal D:
(P_u)_\#\nu_*\not\ll\mathcal L^1\right\}
\leq d_\perp-t<\sigma.
\end{equation}
Here the local parameter chart is bi-Lipschitz, so dimensions
and Frostman bounds transfer between the chart and the cap.
The exceptional set has zero $\widetilde\lambda$-mass.
Since $\lambda$ is carried by the compact set
$L_{\mathbf a}(K)$, we may choose a nonexceptional
$v\in\mathcal D$ and $a_{r-1}\in K$ satisfying
$L_{\mathbf a}(a_{r-1})=\mathbb R v$.

Set $F:=\aff(a_0,\ldots,a_{r-1})=A+\mathbb R v$.
Orthogonal decomposition gives
\begin{equation}\label{eq:height-identity}
P_v(\Pi_A(y))=\text{dist}(\Pi_A(y),\R v)=\dist(y,F).
\end{equation}
Thus the nonzero absolutely continuous measure
$(P_v)_\#\nu_*$ is carried by
$\{\dist(y,F):y\in E\}$, which therefore has positive
Lebesgue measure. Finally,
\begin{equation}\label{eq:base-times-height}
\operatorname{Vol}_r(a_0,\ldots,a_{r-1},y)
=\frac{\operatorname{Vol}_{r-1}(a_0,\ldots,a_{r-1})}{r}
\dist(y,F).
\end{equation}

The coefficient is positive by good position. Hence
\[
\mathcal L^1
\left(
\{\operatorname{Vol}_r(a_0,\ldots,a_{r-1},y):y\in E\}
\right)>0.
\]
Since this set is contained in $V_r(E)$, we also have
$\mathcal L^1(V_r(E))>0$.
\end{proof}

Next, we prove Corollary \ref{cor:almost-every-base} recycling some notation and ideas from the proof above.

\begin{proof}[Proof of Corollary \ref{cor:almost-every-base}]
Fix $s$ as in the statement. By compact approximation,
we may replace $E$ by a compact subset with Hausdorff
dimension greater than $s$. Proving the conclusion for
this subset suffices for the original set.

Choose an $s$-Frostman probability measure $\mu$ on $E$.
As in the proof of Theorem~\ref{thm:main}, restrict $\mu$
to an affine subspace $H$ of minimal dimension carrying
positive mass, and normalize. The resulting measure
is irreducible in $H$ and remains $s$-Frostman.
Since $s>r-1$, we have $\dim H\geq r$.
Identify $H$ isometrically with a Euclidean space,
again denoting its dimension by $d$, and replace $E$
by the compact support of the restricted measure.
The inequality
\[
s>\max\left\{r-1,\frac{d-r+2}{2}\right\}
\]
continues to hold.

Restricting and normalizing $\mu$ in sufficiently small
closed balls around $r$ affinely independent support
points gives irreducible $s$-Frostman probability
measures $\mu_0,\ldots,\mu_{r-1}$ in good position.
Write $E_i=\supp\mu_i$, and let
\[
\mathcal N
:=
\left\{
(a_0,\ldots,a_{r-1})\in E_0\times\cdots\times E_{r-1}:
\mathcal L^1
\left(
\{\operatorname{Vol}_r(a_0,\ldots,a_{r-1},x):x\in E\}
\right)=0
\right\}.
\]

To see that $\mathcal N$ is Borel, write
\[
S_{\mathbf a}
:=
\{\operatorname{Vol}_r(a_0,\ldots,a_{r-1},x):x\in E\}.
\]
If $\mathcal L^1(S_{\mathbf a})<\eta$, compactness and
outer regularity give an open set $O\supset S_{\mathbf a}$
with $\mathcal L^1(O)<\eta$. By uniform continuity of
the volume function, $S_{\mathbf b}\subset O$ for all
bases $\mathbf b$ sufficiently close to $\mathbf a$.
Thus $\{\mathbf a:\mathcal L^1(S_{\mathbf a})<\eta\}$
is relatively open in $E_0\times\cdots\times E_{r-1}$.
Since
\[
\mathcal N
=
\bigcap_{n=1}^{\infty}
\{\mathbf a:\mathcal L^1(S_{\mathbf a})<1/n\},
\]
the set $\mathcal N$ is Borel.

Suppose, towards a contradiction, that
\[
(\mu_0\times\cdots\times\mu_{r-1})(\mathcal N)>0.
\]
Choose
\[
0<\varepsilon
<
(\mu_0\times\cdots\times\mu_{r-1})(\mathcal N).
\]
If $d>r$, choose $\sigma$ as in \eqref{eq:sigma-choice}.
If $d=r$, we are assuming by hypothesis that $\dim_H(E)>d-1$. Choose
\[
d-s<\sigma<1,
\]
which is possible because $s>d-1$.
In both cases,
\[
0<\sigma<\min\{s,d-r+1\},
\]
so Lemma~\ref{lem:boz-thin-input} applies. The resulting
set $G$ satisfies
\[
(\mu_0\times\cdots\times\mu_{r-1})(G\cap\mathcal N)>0.
\]
By Fubini, there is an anchor tuple
$\mathbf a=(a_0,\ldots,a_{r-2})$ for which
\[
\mu_{r-1}\bigl((G\cap\mathcal N)_{\mathbf a}\bigr)>0.
\]
Choose a compact set
\[
K\subset(G\cap\mathcal N)_{\mathbf a},
\qquad \mu_{r-1}(K)>0.
\]
Define $\tilde{\mu}$, $A$, $\nu$, and $\lambda$ from this
anchor and this set $K$ as in the proof of
Theorem~\ref{thm:main}.

If $d>r$, the remainder of that proof applies unchanged.
It produces $a_{r-1}\in K$ for which the corresponding
base-pinned volume set has positive Lebesgue measure,
contradicting $K\subset\mathcal N_{\mathbf a}$.

It remains to consider $d=r$, when the normal space
of the anchor has dimension two. The direction measure
$\lambda$ is $\sigma$-Frostman by
\eqref{eq:direction-thin-frostman}. The projected vertex
measure $\nu$ satisfies the estimate
\[
\nu(B(z,\delta))\lesssim_K\delta^{s-d+2},
\qquad 0<\delta\leq1.
\]
Indeed, the intersection of $K$ with the preimage of
such a ball can be covered by
$O_K(\delta^{-(r-2)})$ balls of radius $O(\delta)$,
and $\tilde{\mu}$ is $s$-Frostman.

Note that $\sigma<1$ while $s-d+2>1$, so $\nu$ has better Frostman parameter than if we had used the thin-plane bound. This is a special feature of the case $s>d-1$. When $s\leq d-1$
observe that $\min\{s,d-r+1\}\geq s-r+2$.

The local transversality argument also applies in
dimension two. Localize $\nu$ and $\lambda$ as in the
main proof, obtaining a probability measure $\nu_*$
and a nonzero direction measure $\widetilde\lambda$.
Since $\sigma>d-s$, we may choose
\[
2-\sigma<t<s-d+2.
\]
In particular, $t>1$, and the
$(s-d+2)$-Frostman bound gives $I_t(\nu_*)<\infty$.
The Peres--Schlag theorem therefore bounds the
exceptional directions by dimension
\[
2-t<\sigma.
\]
They have zero $\widetilde\lambda$-mass, so the same
selection argument produces $a_{r-1}\in K$ for which
the base-pinned volume set has positive Lebesgue measure.
This again contradicts $K\subset\mathcal N_{\mathbf a}$.
\end{proof}

\section{Direct connections with the pinned Falconer problem}
\label{sec:falconer-connection}

The relationship between pinned simplex volumes and pinned distances is more
than motivational. For each fixed base, it is an exact identity after
orthogonal projection to the normal space of the base. There is also an exact
converse on a natural layered class, provided that the base is required to
have a matched form. We record both statements and then give a lattice obstruction to a
universal volume threshold.

\subsection{A fixed base is a pinned distance problem}

\begin{proposition}
\label{prop:normal-space-distance}
Let $2\leq r\leq d$, let $E\subset\mathbb R^d$, and let
$x_0,\ldots,x_{r-1}$ be affinely independent points. Put
\[
U=
\operatorname{span}\{x_1-x_0,\ldots,x_{r-1}-x_0\},
\]
let
\[
\pi_{U^\perp}:\mathbb R^d\longrightarrow U^\perp
\]
be orthogonal projection, and set
\[
z=\pi_{U^\perp}(x_0).
\]
Then
\begin{equation}
\begin{aligned}
&\{\operatorname{Vol}_r(x_0,\ldots,x_{r-1},x):x\in E\}
\\
&\qquad =
\frac{\operatorname{Vol}_{r-1}(x_0,\ldots,x_{r-1})}{r}
\Delta_z(\pi_{U^\perp}(E)).
\end{aligned}
\label{eq:normal-space-distance}
\end{equation}
In particular, the pinned simplex-volume set is a fixed positive dilation of
an ordinary pinned distance set in the $(d-r+1)$-dimensional space
$U^\perp$.
\end{proposition}

Before proving the proposition above, note that the identity \eqref{eq:normal-space-distance} suggests a natural conjectural threshold for pinned simplex volumes. Indeed, the pinned Falconer conjecture suggests $\frac{d-r+1}{2}$ as the relevant dimensional threshold in the normal space, while the condition $\dim_{\mathrm H}(E)>r-1$ is necessary to exclude sets contained in affine $(r-1)$-planes. This motivates the following conjecture.

\begin{conjecture}
\label{conj:pinned-simplex-volumes}
Let $2\leq r\leq d$ and let $E\subset\mathbb R^d$ be a Borel set satisfying
\[
\dim_{\mathrm H}(E)>
\max\left\{r-1,\frac{d-r+1}{2}\right\}.
\]
Then there exist affinely independent points $x_0,\ldots,x_{r-1}\in E$ such that
\[
\mathcal L^1\left(
\left\{\operatorname{Vol}_r(x_0,\ldots,x_{r-1},x):x\in E\right\}
\right)>0.
\]
\end{conjecture}

Despite the exact correspondence in \eqref{eq:normal-space-distance}, a pinned distance theorem cannot be applied directly: both the projected set $\pi_{U^\perp}(E)$ and the prescribed pin $z=\pi_{U^\perp}(x_0)$ depend on the chosen base. In particular, even if the projected set has sufficient dimension, $z$ is not necessarily a good pin. Our approach instead exploits the geometry of affine planes determined by $E$ to select a base with a positive-measure height set. Theorem~\ref{thm:main} establishes the conjectured conclusion under the stronger hypothesis
\[
\dim_{\mathrm H}(E)>
\max\left\{r-1,\frac{d-r+2}{2}\right\}.
\]

\begin{proof}[Proof of Proposition \ref{prop:normal-space-distance}]
Let
\[
F=x_0+U
\]
be the affine span of the base. Since $\pi_{U^\perp}$ annihilates $U$, every point of
$F$ has the same image $z$. Orthogonal decomposition gives
\[
\operatorname{dist}(x,F)=|\pi_{U^\perp}(x)-z|.
\]
The base-times-height identity gives
\[
\operatorname{Vol}_r(x_0,\ldots,x_{r-1},x)
=
\frac{\operatorname{Vol}_{r-1}(x_0,\ldots,x_{r-1})}{r}
\operatorname{dist}(x,F).
\]
Combining the two identities and allowing $x$ to vary in $E$ proves
\eqref{eq:normal-space-distance}.
\end{proof}

The same reduction can be written in graph coordinates. Write
\[
\mathbb R^d=\mathbb R^{r-1}\times\mathbb R^c,
\qquad
c=d-r+1,
\]
and write a nearby affine base plane as
\[
L_{A,b}=\{(p,Ap+b):p\in\mathbb R^{r-1}\},
\]
where $A$ is a real $c$ by $(r-1)$ matrix and $b\in\mathbb R^c$. Let
\[
M_A=(I_c+AA^T)^{-\frac12},
\]
where $I_c$ is the $c$ by $c$ identity matrix and the right side is the
positive-definite inverse square root. Define
\[
T_A(p,q)=M_A(q-Ap)
\]
and
\[
z_{A,b}=M_Ab.
\]
The rows of $T_A$ are orthonormal, since
$M_A(I_c+AA^T)M_A=I_c$, and its kernel is the direction space
of $L_{A,b}$. Consequently,
\begin{equation}
\operatorname{dist}((p,q),L_{A,b})
=
|T_A(p,q)-z_{A,b}|.
\label{eq:fixed-orientation-distance}
\end{equation}
Thus, for each fixed orientation $A$, translation of the affine base plane is
exactly translation of an ordinary Falconer pin in $\mathbb R^c$. If
$x_0,\ldots,x_{r-1}$ span $L_{A,b}$, then
\[
T_A(x_0)=\cdots=T_A(x_{r-1})=z_{A,b}.
\]

\subsection{Layered lifts and matched bases}

We next show that every pinned distance set is itself a pinned simplex-volume
set of a special form. Let $n\geq1$ and $r\geq2$. In $\mathbb R^{r-1}$, let
\[
e_0=0
\]
and let $e_1,\ldots,e_{r-1}$ be the standard basis vectors. For
$E\subset\mathbb R^n$, define the layered lift
\begin{equation}
\widetilde E_r
=
\bigcup_{j=0}^{r-1}(E\times\{e_j\})
\subset\mathbb R^{n+r-1}.
\label{eq:layered-lift}
\end{equation}
For each $x\in E$, define the matched base
\begin{equation}
B_x=
\bigl((x,e_0),(x,e_1),\ldots,(x,e_{r-1})\bigr).
\label{eq:matched-base}
\end{equation}

\begin{proposition}
\label{prop:layered-lift}
Let $E\subset\mathbb R^n$. Then
\[
\dim_{\mathrm H}(\widetilde E_r)=\dim_{\mathrm H}(E).
\]
Moreover, for every $x\in E$,
\begin{equation}
\begin{aligned}
&\{\operatorname{Vol}_r((x,e_0),\ldots,(x,e_{r-1}),z):
 z\in\widetilde E_r\}
\\
&\qquad =\frac{1}{r!}\Delta_x(E).
\end{aligned}
\label{eq:layered-distance-identity}
\end{equation}
\end{proposition}

\begin{proof}
The set $\widetilde E_r$ is a finite union of isometric copies of $E$, so it
has the same Hausdorff dimension as $E$.

Every $z\in\widetilde E_r$ has the form
\[
z=(y,e_j)
\]
for some $y\in E$ and some $0\leq j\leq r-1$. Using $(x,e_0)$ as the initial
vertex, the edge vectors of the corresponding $r$-simplex are
\[
(0,e_1),\ldots,(0,e_{r-1}),(y-x,e_j).
\]
The vertical component $(0,e_j)$ of the last vector lies in the span of the
first $r-1$ vectors and therefore makes no contribution to their exterior
product. Here $\wedge$ denotes the exterior product, and $|\cdot|$ denotes
its Euclidean norm. Hence
\[
\begin{aligned}
&r!\operatorname{Vol}_r((x,e_0),\ldots,(x,e_{r-1}),(y,e_j))
\\
&\qquad =
\left|
(0,e_1)\wedge\cdots\wedge(0,e_{r-1})\wedge(y-x,0)
\right|.
\end{aligned}
\]
The first $r-1$ vectors are orthonormal and are perpendicular to $(y-x,0)$.
The exterior product therefore has norm $|y-x|$, and
\[
\operatorname{Vol}_r((x,e_0),\ldots,(x,e_{r-1}),(y,e_j))
=
\frac{|y-x|}{r!}.
\]
Allowing $y$ and $j$ to vary proves \eqref{eq:layered-distance-identity}.
\end{proof}

Proposition \ref{prop:layered-lift} gives an exact equivalence at the level of
sharp conjectures. Fix $n\geq2$ and $r\geq2$. The assertion that every Borel
set $E\subset\mathbb R^n$ with
\[
\dim_{\mathrm H}(E)>\frac n2
\]
has a point $x\in E$ for which $\Delta_x(E)$ has positive Lebesgue measure is
equivalent to the assertion that every layered lift $\widetilde E_r$ of such
a set has a matched base $B_x$ for which the pinned $r$-simplex volume set has
positive Lebesgue measure. If
\[
D=n+r-1
\]
is the ambient dimension of the lift, then the required threshold on this
special class is
\begin{equation}
\frac n2=\frac{D-r+1}{2}.
\label{eq:layered-sharp-threshold}
\end{equation}
There are two obstacles to deducing the sharp pinned Falconer
statement from Theorem~\ref{thm:main}: the base must have the
matched form, and the dimensional threshold must be lower. Indeed, applying Theorem~\ref{thm:main} to the lifted set
$\widetilde E_r\subset\mathbb R^{n+r-1}$ gives a good
base under the hypothesis
$$\dimH E>
\max\left\{r-1,\frac{D-r+2}{2}\right\}
=
\max\left\{r-1,\frac{n+1}{2}\right\}.$$
This threshold is strictly greater than $n/2$.
Thus, even if the theorem selected a matched base,
its dimensional hypothesis would not yield the
conjectured pinned Falconer threshold.

The matched-base requirement is essential. Theorem~\ref{thm:main}
produces some good base, but when applied to $\widetilde E_r$,
there is no reason for the selected base to
be one of the bases $B_x$. Already for triangles, take
\[
a=(u,0),
\qquad
b=(v,1),
\qquad
z=(y,0)
\]
in $\mathbb R^{n+1}$. The Gram determinant gives
\begin{equation}
\begin{aligned}
4\operatorname{Vol}_2(a,b,z)^2
&=
(1+|v-u|^2)|y-u|^2
\\
&\qquad -((v-u)\cdot(y-u))^2.
\end{aligned}
\label{eq:unmatched-base-anisotropic}
\end{equation}
Unless $v=u$, this is an anisotropic quadratic norm of $y-u$, rather than the
Euclidean distance from $y$ to $u$. Thus positive measure for the volume set
associated with an arbitrary cross-layer base does not give the pinned
distance conclusion in \eqref{eq:layered-distance-identity}.

The same point remains even for an abundance theorem about arbitrary bases.
If $\mu$ is a non-atomic measure on $E$, the matched bases correspond to the
diagonal
\[
\{(x,\ldots,x):x\in E\}
\]
in the $r$ horizontal coordinates. This diagonal has measure zero with respect to the $r$-fold product measure
$\mu\times\cdots\times\mu$. Consequently, a statement holding for almost
every base with respect to a product Frostman measure need not produce a
single matched base. What would
suffice is an exceptional-base theorem showing that the matched family cannot
be entirely exceptional.

\subsection{A Falconer-type lattice obstruction}

The classical lattice obstruction for the Falconer distance problem also has
a direct simplex-volume analogue. We first record the rapidly nested lattice
construction in a form that will be used in the proof.

\begin{lemma}
\label{lem:rapid-lattice-set}
Let $m\geq1$ and let $\alpha>1$. There are a compact set
$F\subset[0,1]^m$, a constant $C_m>0$, and integers $q_j$ tending to infinity
such that
\[
\dim_{\mathrm H}(F)=\frac{m}{\alpha}
\]
and, for every $j$,
\begin{equation}
F\subset
\left\{x\in[0,1]^m:
\operatorname{dist}(x,q_j^{-1}\mathbb Z^m)
\leq C_mq_j^{-\alpha}
\right\}.
\label{eq:rapid-lattice-approximation}
\end{equation}
\end{lemma}

\begin{proof}
Choose the positive integers $q_j$ recursively so rapidly that
\[
q_jq_{j-1}^{-\alpha}>100
\]
for $j\geq2$ and
\begin{equation}
\frac{\sum_{i=1}^{j-1}\log q_i}{\log q_j}
\longrightarrow0.
\label{eq:rapid-growth}
\end{equation}
We may also require that $q_j$ is a multiple of $q_{j-1}$. Put
\[
\ell_j=q_j^{-\alpha}.
\]
Since $\alpha>1$, after increasing $q_1$ if necessary, one has
\[
2\ell_j<\frac{1}{4q_j}
\]
for every $j$.

At the first stage, take a collection of pairwise disjoint closed cubes of
side length $2\ell_1$, centered at points of $q_1^{-1}\mathbb Z^m$, and
contained in $[0,1]^m$. At stage $j\geq2$, inside each cube from stage $j-1$
choose the same number $M_j$ of pairwise disjoint closed cubes of side length
$2\ell_j$, centered at points of $q_j^{-1}\mathbb Z^m$, and contained in the
parent cube. The rapid growth of $q_j$ allows this number to be chosen so that
there are constants $c_m,C_m'>0$, independent of $j$, for which
\begin{equation}
c_m(q_j\ell_{j-1})^m
\leq M_j\leq
C_m'(q_j\ell_{j-1})^m.
\label{eq:child-count}
\end{equation}
At the first stage, choose $M_1$ cubes with
\[
c_mq_1^m\leq M_1\leq C_m'q_1^m.
\]
Let $F_j$ be the union of all stage $j$ cubes and put
\[
F=\bigcap_{j=1}^{\infty}F_j.
\]
The sets $F_j$ are nonempty compact sets and are nested, so $F$ is nonempty
and compact. Every point of $F$ lies in a stage $j$ cube centered at a point
of $q_j^{-1}\mathbb Z^m$. This proves
\eqref{eq:rapid-lattice-approximation}, after changing the dimensional
constant.

Let
\[
N_j=M_1\cdots M_j
\]
be the number of stage $j$ cubes. Taking logarithms in
\eqref{eq:child-count} gives a constant $C>0$ such that
\[
\left|
\log N_j-
\left(
m\log q_j-(\alpha-1)m\sum_{i=1}^{j-1}\log q_i
\right)
\right|
\leq Cj.
\]
It follows from \eqref{eq:rapid-growth} that
\begin{equation}
\frac{\log N_j}{-\log\ell_j}
\longrightarrow
\frac{m}{\alpha}.
\label{eq:moran-dimension-ratio}
\end{equation}
Covering $F$ by its $N_j$ stage $j$ cubes gives
\[
\dim_{\mathrm H}(F)\leq\frac{m}{\alpha}.
\]

For the reverse inequality, assign mass $N_j^{-1}$ to every stage $j$ cube.
The consistency of the equal branching construction defines a Borel
probability measure $\nu$ supported on $F$. Fix
\[
0<\beta<\frac{m}{\alpha}.
\]
Equation \eqref{eq:moran-dimension-ratio} gives, for all sufficiently large
$j$,
\begin{equation}
N_j^{-1}\leq\ell_j^\beta.
\label{eq:level-mass-bound}
\end{equation}
Let $\ell_j\leq\rho<\ell_{j-1}$. If
\[
\rho<\frac{1}{4q_j},
\]
then a ball of radius $\rho$ meets at most a bounded number of stage $j$
cubes, because their centers lie on the $q_j^{-1}$ lattice. Hence
\[
\nu(B(x,\rho))
\leq
C_mN_j^{-1}
\leq
C_m\rho^\beta.
\]
If instead
\[
\frac{1}{4q_j}\leq\rho<\ell_{j-1},
\]
then the ball meets at most a bounded number of stage $j-1$ cubes and at most
\[
C_m(\rho q_j)^m
\]
stage $j$ cubes inside each one. By \eqref{eq:child-count},
\[
\begin{aligned}
\nu(B(x,\rho))
&\leq
C_mN_{j-1}^{-1}\frac{(\rho q_j)^m}{M_j}
\\
&\leq
C_mN_{j-1}^{-1}
\left(\frac{\rho}{\ell_{j-1}}\right)^m.
\end{aligned}
\]
Since $\beta<m$ and $\rho<\ell_{j-1}$, this is at most
\[
C_mN_{j-1}^{-1}
\left(\frac{\rho}{\ell_{j-1}}\right)^\beta
\leq
C_m\rho^\beta
\]
by \eqref{eq:level-mass-bound}. Enlarging the constant handles the finitely
many remaining scales. The mass distribution principle gives
\[
\dim_{\mathrm H}(F)\geq\beta.
\]
Letting $\beta$ increase to $\frac{m}{\alpha}$ completes the proof.
\end{proof}

\begin{proposition}
\label{prop:lattice-obstruction}
Let $1\leq r\leq d$. For every
\[
0<s<\frac{d}{2r},
\]
there is a compact set $E\subset[0,1]^d$ such that
\[
\dim_{\mathrm H}(E)=s
\]
and
\[
\mathcal L^1(V_r(E))=0.
\]
Consequently, any universal Hausdorff-dimension threshold forcing a
positive-measure unpinned or pinned $r$-simplex volume set must be at least
\begin{equation}
\max\left\{r-1,\frac{d}{2r}\right\}.
\label{eq:lattice-lower-bound}
\end{equation}
For $r=1$, this is Falconer's classical lattice obstruction
$\frac d2$.
\end{proposition}

\begin{proof}
Choose
\[
\alpha=\frac ds>2r
\]
and apply Lemma \ref{lem:rapid-lattice-set} with $m=d$. Let $E$ be the
resulting compact set and let $q=q_j$ be one of the corresponding lattice
scales. Every point of $E$ lies within $C_dq^{-\alpha}$ of a point of
$q^{-1}\mathbb Z^d$ in a fixed bounded region.

Consider first an exact lattice simplex with vertices
\[
y_0,\ldots,y_r\in q^{-1}\mathbb Z^d
\]
in that region. Put
\[
v_i=q(y_i-y_0)\in\mathbb Z^d,
\qquad
1\leq i\leq r.
\]
The Gram determinant formula gives
\begin{equation}
(r!)^2\operatorname{Vol}_r(y_0,\ldots,y_r)^2
=
q^{-2r}\det\bigl[(v_i\cdot v_j)_{i,j=1}^r\bigr].
\label{eq:lattice-gram-determinant}
\end{equation}
The determinant on the right is a nonnegative integer. Since the vectors
$v_i$ have length at most $C_dq$, Hadamard's inequality bounds it by
$C_{d,r}q^{2r}$. It follows that the exact lattice determines at most
\[
C_{d,r}q^{2r}
\]
distinct $r$-simplex volumes.

The volume map is Lipschitz on bounded subsets of
$(\mathbb R^d)^{r+1}$. Indeed,
\[
r!\operatorname{Vol}_r(x_0,\ldots,x_r)
=
\left|
(x_1-x_0)\wedge\cdots\wedge(x_r-x_0)
\right|,
\]
and multilinearity followed by a telescoping expansion gives the required
bound. Consequently, every $r$-simplex with vertices in $E$ has volume within
$C_{d,r}q^{-\alpha}$ of one of the exact lattice volumes. Therefore
$V_r(E)$ is covered by at most $C_{d,r}q^{2r}$ intervals of length
$C_{d,r}q^{-\alpha}$, and
\[
\mathcal L^1(V_r(E))
\leq
C_{d,r}q^{2r-\alpha}.
\]
Letting $q=q_j$ tend to infinity gives
\[
\mathcal L^1(V_r(E))=0.
\]
Lemma \ref{lem:rapid-lattice-set} gives
\[
\dim_{\mathrm H}(E)=\frac d\alpha=s.
\]

The second obstruction in \eqref{eq:lattice-lower-bound} follows from this
construction. The first follows by taking a set of Hausdorff dimension
$r-1$ contained in an affine $(r-1)$-plane, since every $r$-simplex with
vertices in such a set has volume zero. This completes the proof.
\end{proof}

\section*{Disclosure of AI assistance}

The project began with the authors' observation that Ren's radial
projection theorem could lead to the threshold $d/2$ for
base-pinned triangle areas. The key was to verify the appropriate
transversality conditions for a local family of distance-to-line
maps and then apply the Peres-Schlag exceptional-set theorem.
An initial approach to higher-dimensional simplices used the
unpinned continuum Beck theorem of Bright, Ortiz, and Zakharov.
This led to a larger parameter family in the Peres--Schlag
argument and consequently to a weaker dimensional threshold for $r\geq 3$.

In developing the
higher-dimensional arguments, discussions with ChatGPT (GPT-6)
helped explore consequences of the geometric estimates of Bright,
Ortiz, and Zakharov, refine the use of projection and exceptional-set
theorems, and develop arguments leading to progressively improved
dimensional thresholds. These discussions contributed to the
development and simplification of the final proof.
ChatGPT also assisted with polishing the manuscript.
The authors take full responsibility for the mathematical
arguments and the final manuscript.

\end{document}